\documentclass[12pt]{article}

\usepackage{tikz}
\usepackage{pgf}

\newtheorem{thm}{Theorem}[section]
\newtheorem{clm}[thm]{Claim}
\newtheorem{cnj}[thm]{Conjecture}
\newtheorem{cor}[thm]{Corollary}

\newtheorem{lem}[thm]{Lemma}

\newtheorem{prp}[thm]{Proposition}
\newtheorem{qst}[thm]{Question}

\usepackage{graphicx,latexsym,amsmath,amsfonts,amssymb,mathrsfs}
\usepackage{ulem} 
\usepackage{xcolor} 

\newcommand{\qbin}[2]{\mathbf{\begin{bmatrix}#1\\#2\end{bmatrix}}}
\newcommand{\qmult}[5]{\mathbf{\begin{bmatrix}&&#1&&\\#2&#3&\ldots&#4&#5\end{bmatrix}}}
\def\BlS{B_{\ell,S}}
\newcommand{\omitt}[1]{}

\def\bid{{\mathbf{id}}}

\def\bG{{\mathbf G}}

\def\cA{{\cal A}}
\def\cB{{\cal B}}
\def\cC{{\cal C}}
\def\cE{{\cal E}}
\def\cF{{\cal F}}  

\def\cG{{\cal G}}
\def\cH{{\cal H}}

\def\cL{{\cal L}}

\def\cP{{\cal P}}

\def\oB{{\overline{B}}}

\def\op{{\overline{p}}}
\def\ox{{\overline{x}}}

\def\sB{{\mathscr B}}

\def\sL{{\mathscr L}}

\def\sP{{\mathscr P}}

\def\Bn{{\sB(n)}}
\def\Bfour{{\sB(4)}}
\def\Bln{{\sB_{\ell}(n)}}
\def\Btwosix{{\sB_2(6)}}
\def\Bthreesix{{\sB_3(6)}}
\def\Bkn{{\sB_k(n)}}
\def\Brn{{\sB_r(n)}}
\def\Bonen{{\sB_1(n)}}
\def\Btwon{{\sB_2(n)}}
\def\Bthreen{{\sB_3(n)}}

\def\cHmrk{{\cH_{m,r,k}}}

\def\dn{{\sf down}}

\def\min{{\sf min}}
\def\rank{{\sf rank}}
\def\rev{{\sf ID}} 

\def\ID{{\sf ID}}
\def\Inv{{\sf Inv}}
\def\up{{\sf up}}

\def\symn{{\rm Sym}(n)}

\newcommand{\precdot}{\prec\mathrel{\mkern-5mu}\mathrel{\cdot}}
\def\mt{{\emptyset}}
\def\rel{{\ \preceq\ }}
\def\crel{{\ \precdot\ }} 
\def\sse{{\ \subseteq\ }}

\def\pr{{\prime}}

\def\multiset#1#2{\ensuremath{\left(\kern-.3em\left(\genfrac{}{}{0pt}{}{#1}{#2}\right)\kern-.3em\right)}}

\def\sqr#1#2{{\vcenter{\hrule height.#2pt
        \hbox{\vrule width.#2pt height#1pt \kern#1pt
                \vrule width.#2pt}
        \hrule height.#2pt}}}

\def\pf{{\hfill$\Box$\\ \medskip}}
\def\proof{{\noindent {\it Proof}.\ \ }}

\providecommand{\binom}[2]{{#1\choose#2}}

\begin{document}

%
%
\title{Erd\H{o}s-Ko-Rado theorems on the weak Bruhat lattice}

\author{
Susanna Fishel\footnote{School of Mathematical and Statistical
  Sciences, Arizona State University, Tempe, AZ,  \tt{sfishel1@asu.edu}.},
\and Glenn Hurlbert\footnote{Virginia Commonwealth University, Richmond, Virginia, \tt{ghurlbert@vcu.edu}.}, 
\and Vikram Kamat\footnote{Villanova University, Villanova, Pennsylvania, \tt{vikram.kamat@villanova.edu}.}, 
\and Karen Meagher \footnote{University of Regina, Regina, Saskatchewan, \tt{karen.meagher@uregina.ca}.}
}


\maketitle

\vspace{1.0 in}

%
%
\begin{abstract}
  Let $\sL=(X,\rel)$ be a lattice. For $\cP\sse X$ we say that $\cP$ is $t$-{\it intersecting} if $\rank(x\wedge y)\ge t$ for all $x,y\in\cP$. The seminal theorem of Erd\H{o}s, Ko and Rado describes the maximum intersecting $\cP$ in the lattice of subsets of a finite set with the additional condition that $\cP$ is contained within a level of the lattice. The Erd\H{o}s-Ko-Rado theorem has been extensively studied and generalized to other objects and lattices.

 In this paper, we focus on intersecting families of permutations as defined with respect to the weak Bruhat lattice. In this setting, we prove analogs of certain extremal results on intersecting set systems. In particular we give a characterization of the maximum intersecting families of permutations in the Bruhat lattice. We also characterize the maximum intersecting families of permutations within the $r^{\textrm{th}}$ level of the Bruhat lattice of permutations of size $n$, provided that $n$ is large relative to $r$.
\end{abstract}


\newpage

%
%
\section{Introduction}\label{sec:intro}

Let $\sL=(X,\rel)$ be a lattice.  For $\cP\sse X$ we say that $\cP$ is
$t$-{\it intersecting} if $\rank(x\wedge y)\ge t$ for all $x,y\in\cP$.
If $\sL$ is the subset lattice, two subsets $A,B$ are
$t$-intersecting exactly when $|A \cap B | \ge t$ (we will refer to $1$-intersecting sets simply as intersecting sets). 
The problem of finding maximum collections of intersecting sets has a long
history, see~\cite{EKRbook} and the references within. Perhaps the most famous result is the
Erd\H{o}s-Ko-Rado (EKR) Theorem~\cite{EKR}. This theorem gives the size of the largest
collection of sets at level $r$ in the subset lattice such that any
two of the sets intersect.

\begin{thm}\label{thm:ekr}
	Let $r$ and $n$ be integers with $n \geq 2r$.
	If $\cA$ is an intersecting collection of $r$-subsets of the set $\{1,\dots,n\}$
	then
	\[
		|\cA| \le \binom{n-1}{r-1}.
	\]
	Moreover, if $n>2r$ equality holds if and only if $\cA$
        consists of all the $r$-subsets that contain a common element.
\end{thm}  

Another well-known result is Katona's theorem~\cite{MR0168468}, which
gives the size of the largest collection of $t$-intersecting sets
(with no restriction on the size of the sets) in the subset lattice on
$n$ points.  If $n+t = 2v$, then this maximum size is achieved by the
collection of all subsets with size larger than $v$. If $n+t = 2v-1$,
an example of a collection of maximum size is the collection that
contains all subsets with size at least $v$, along with all subsets of
size $v-1$ that do not contain a fixed element. In the special case
that $t=1$ these sets have size $2^{n-1}$. This implies that the
collection of all subsets that contain a fixed element is a
$1$-intersecting set of maximum size. In this case there are also many
other (non-isomorphic) collections of maximum size.

One of the reasons that the EKR Theorem is so famous is that versions
of it hold for many other objects; for example there are versions for
integer sequences, vector spaces over a finite field, perfect
matchings, independent sets in graphs,
and permutations (these are just a few of the examples;
again see~\cite{EKRbook}, and the references within, for more details).

In this paper, we give versions of Katona's Theorem and the
Erd\H{o}s-Ko-Rado Theorem for the weak Bruhat lattice of
permutations. We start with some general results for all
lattices. Next we define and give properties of the weak Bruhat
lattice. Section 4 has information about properties of intersecting
sets in the Bruhat lattice and a characterization for the maximum
intersecting collections of intersecting permutations in the Bruhat
lattice. Section 5 gives a version of the EKR theorem for permutations
in level $r$ of the Bruhat lattice, provided that $n$ is large
relative to $r$.  In the sequel, collections of sets will be called
{\it systems}, and collections of permutations will be called {\it
  families}.

\section{General lattices}\label{Defs}

For a subset $\cF$ of elements in a lattice $\sL=(X,\rel)$, define the {\it upset
  of $\cF$} to be
\[
\up(\cF)=\{z\in X\mid x\rel z, {\ \rm some\ } x\in\cF\}.
\]
Similarly, the {\it downset of $\cF$} is defined to be
\[
\dn(\cF)=\{z\in X\mid z\rel x, {\ \rm some\ } x\in\cF\}.
\]
For any $\cF$, if $\up(\cF)=\cF$, then we say that $\cF$ is an
{\it upset}. Similarly, $\cF$ is a {\it downset} if
$\dn(\cF)=\cF$. Upsets and downsets are also called filters and
ideals, respectively.

For a single element $p\in \sL$, the family $\up(\{p\})$ (which we
abbreviate $\up(p)$) is called a {\it star} with {\it center $p$}. If
the rank of $p$ is $t$, then $\up(p)$ is $t$-intersecting; we call
$\up(p)$ a {\it canonical} $t$-intersecting family. More generally,
for some $\cH \sse \sL$ we define $\up_\cH(\cF)=\up(\cF)\cap\cH$.  In
this paper, we only consider the case where $\cH$ is the set of all
elements in $\sL$ with rank $k$; this set is denoted by $\sL_k$.  If
$\rank(p) =t$, then $\up_{\sL_k }(p)$ is called a {\it $t$-star at
  level $k$}.

Define $f_t(\sL)$ to be the maximum size of a $t$-intersecting family
from the lattice $\sL$. We say that $\sL$ has the {\it $t$-EKR
  property} (or is {\it $t$-EKR}) if $f_t(\sL)$ is equal to the size of a
$t$-star. Equivalently, $\sL$ has the $t$-EKR property if $f_t(\sL) =
\max_{p \in \sL_t} |\up(p)|$. A lattice that has the $1$-EKR property
will simply be said to have the EKR property.

Similarly, $f_t(\sL_k)$ is defined to be the maximum size of a
$t$-intersecting family in level $k$ of the lattice. The
level $\sL_k$ has the $t$-{\it EKR property} (or {\it is $t$-EKR}) if
the size of the largest $t$-intersecting family of $\sL_k$ can be
realized by a $t$-star at level $k$. This is equivalent to saying that
$f_t(\sL_k) = \max_{p \in \sL_t} |\up_{\cL_k}(p)|$.
Again, we suppress the use of $t$ when $t=1$.

The Erd\H{o}s-Ko-Rado Theorem is equivalent to saying that in the
subset lattice all levels below $n/2$ have the $1$-EKR
property. Wilson~\cite{MR0771733} proved the subset lattice for
$\{1,\ldots,n\}$ has the $t$-EKR property at level-$k$ provided that
$n>(t+1)(k-t+1)$ (Frankl~\cite{MR519277} had previously proved this
with the restriction that $t>14$). Ahlswede and Khachatrian's complete
intersection theorem~\cite{MR1429238} describes all the largest
intersecting sets for all values of $n$, $k$ and $t$.

We say that a lattice $\sL$ is {\it uniquely complemented} if, for all
$x\in X$, there exists a unique $y\in X$ such that
$\rank(x\wedge y)=0$ and $\rank(x\vee y)=\rank(\sL)$. 
A lattice is {\it distributive} if for all $x,y,z \in \sL$ the
following holds
\[
x \wedge (y \vee z) = (x \wedge y) \vee (x \wedge z).
\]

\begin{thm}\label{thm:complemented}
Any uniquely complemented distributive lattice has the EKR property.
\end{thm}

\proof
Assume that $\sL$ is a uniquely complemented lattice and denote the
least element in $\sL$ by $0$. 
Denoting the complement of an element $x\in\sL$ by $\ox$,
at most one of the pair $\{x,\ox\}$ can belong to an
intersecting family. Thus an intersecting set can be no larger than $|\sL|/2$.

Let $p$ be any element of the lattice with rank 1, and assume that for
some $x \in \sL$, both $x$ and $\ox$ are not comparable to $p$. Then
\[
0 = (p \wedge x) \vee (p \wedge \ox) = p \wedge (x \vee \ox) = p
\]                      
which is a contradiction. Thus for any element $x \in \sL$, either $x$ or $\ox$ is comparable to $p$. Since $p$ is rank 1 this means that either $x$ or $\ox$ is above $p$ in the lattice, and so
$|\up(p)| = |\sL|/2$.
\pf

Since the subset lattice is uniquely complemented, this theorem
implies that the subset lattice has the EKR property.  Katona's
theorem shows that the subset lattice does not have the $t$-EKR
property for $t >1$. 

It seems to be harder to find general results about maximum
intersecting families within a level of a lattice.  Many of the standard
proofs for the EKR property use regularity conditions on the
set. Suda~\cite{MR2901181} considered the $t$-EKR property for levels
in semi-lattices with regularity conditions. The Bruhat lattice is
different since it does not satisfy these regularity
conditions. Specifically the upsets of different elements at the same
level may have different sizes, this makes the Bruhat lattice a
particularly interesting lattice to consider. More details are given
in the next section.

\section{Properties of the weak Bruhat lattice}

\begin{figure}
\begin{center}
\begin{tikzpicture}
\def\x{1.25} \def\ye{.25}
\begin{scope}[shift={(0*\x,0)}]
\node (e) at (0,0) {$1234$};
\node (s1) at (-\x,1) {$2134$};
\node (s2) at (0,1) {$1324$};
\node (s3) at (\x,1) {$1243$};

\node (s12) at (-2*\x,2) {$2314$};
\node (s21) at (-\x,2) {$3124$};
\node (s23) at (0,2) {$1342$};
\node (s13) at (\x,2) {$2143$};
\node (s32) at (2*\x,2) {$1423$};

\node (s123) at (-2.5*\x,3) {$2341$};
\node (s121) at (-1.5*\x,3) {$3214$};
\node (s231) at (-.5*\x,3) {$3142$};
\node (s132) at (.5*\x,3) {$2413$};
\node (s232) at (1.5*\x,3) {$1432$};
\node (s321) at (2.5*\x,3) {$4123$};

\node (s1231) at (-2*\x,4) {$3241$};
\node (s1232) at (-\x,4) {$2431$};
\node (s2312) at (0,4) {$3412$};
\node (s1321) at (\x,4) {$4213$};
\node (s2321) at (2*\x,4) {$4132$};

\node (s12312) at (-1.5*\x,5) {$3421$};
\node (s12321) at (0,5) {$4231$};
\node (s23121) at (1.5*\x,5) {$4312$};

\node (s123121) at (0,6) {$4321$};

\draw (e)--(s1);
\draw (e)--(s2);
\draw (e)--(s3);
\draw (s1)--(s12);
\draw (s1)--(s13);
\draw (s2)--(s21);
\draw (s2)--(s23);
\draw (s3)--(s13);
\draw (s3)--(s32);
\draw (s12)--(s123);
\draw (s12)--(s121);
\draw (s21)--(s121);
\draw (s21)--(s231);
\draw (s13)--(s132);
\draw (s23)--(s231);
\draw (s23)--(s232);
\draw (s32)--(s232);
\draw (s32)--(s321);
\draw (s123)--(s1231);
\draw (s123)--(s1232);
\draw (s121)--(s1231);
\draw (s231)--(s2312);
\draw (s132)--(s1232);
\draw (s132)--(s1321);
\draw (s232)--(s2321);
\draw (s321)--(s2321);
\draw (s321)--(s1321);

\draw (s1231)--(s12312);
\draw (s2312)--(s12312);
\draw (s1232)--(s12321);
\draw (s1321)--(s12321);
\draw (s2312)--(s23121);
\draw (s2321)--(s23121);

\draw (s12312)--(s123121);
\draw (s12321)--(s123121);
\draw (s23121)--(s123121);
\end{scope}

\begin{scope}[shift={(6*\x,0)}]
\node (e) at (0,0) {$\{\}$};
\node (s1) at (-\x,1) {$\{1\}$};
\node (s2) at (0,1) {$\{2\}$};
\node (s3) at (\x,1) {$\{3\}$};

\node (s12) at (-2*\x,2) {$\{1\}$};
\node (s21) at (-\x,2) {$\{2\}$};
\node (s23) at (0,2) {$\{2\}$};
\node (s13) at (\x,2) {$\{1,3\}$};
\node (s32) at (2*\x,2) {$\{3\}$};

\node (s123) at (-2.5*\x,3) {$\{1\}$};
\node (s121) at (-1.5*\x,3) {$\{1,2\}$};
\node (s231) at (-.5*\x,3) {$\{2\}$};
\node (s132) at (.5*\x,3) {$\{1,3\}$};
\node (s232) at (1.5*\x,3) {$\{2,3\}$};
\node (s321) at (2.5*\x,3) {$\{3\}$};

\node (s1231) at (-2*\x,4) {$\{1,2\}$};
\node (s1232) at (-\x,4) {$\{1,3\}$};
\node (s2312) at (0,4) {$\{2\}$};
\node (s1321) at (\x,4) {$\{1,3\}$};
\node (s2321) at (2*\x,4) {$\{2,3\}$};

\node (s12312) at (-1.5*\x,5) {$\{1,2\}$};
\node (s12321) at (0,5) {$\{1,3\}$};
\node (s23121) at (1.5*\x,5) {$\{2,3\}$};

\node (s123121) at (0,6) {$\{1,2,3\}$};

\draw (e)--(s1);
\draw (e)--(s2);
\draw (e)--(s3);
\draw (s1)--(s12);
\draw (s1)--(s13);
\draw (s2)--(s21);
\draw (s2)--(s23);
\draw (s3)--(s13);
\draw (s3)--(s32);
\draw (s12)--(s123);
\draw (s12)--(s121);
\draw (s21)--(s121);
\draw (s21)--(s231);
\draw (s13)--(s132);
\draw (s23)--(s231);
\draw (s23)--(s232);
\draw (s32)--(s232);
\draw (s32)--(s321);
\draw (s123)--(s1231);
\draw (s123)--(s1232);
\draw (s121)--(s1231);
\draw (s231)--(s2312);
\draw (s132)--(s1232);
\draw (s132)--(s1321);
\draw (s232)--(s2321);
\draw (s321)--(s2321);
\draw (s321)--(s1321);

\draw (s1231)--(s12312);
\draw (s2312)--(s12312);
\draw (s1232)--(s12321);
\draw (s1321)--(s12321);
\draw (s2312)--(s23121);
\draw (s2321)--(s23121);

\draw (s12312)--(s123121);
\draw (s12321)--(s123121);
\draw (s23121)--(s123121);
\end{scope}

\end{tikzpicture}
\caption{Two representations of the Bruhat lattice $\Bfour$. On the left, vertices are labeled by $p=p_1p_2p_3p_4$, 
where $p(i)=p_i$.  On the right, each vertex $p$ is labeled by its set of inverse descents $\rev(p)$. }
\label{figS4}
\end{center}
\end{figure}
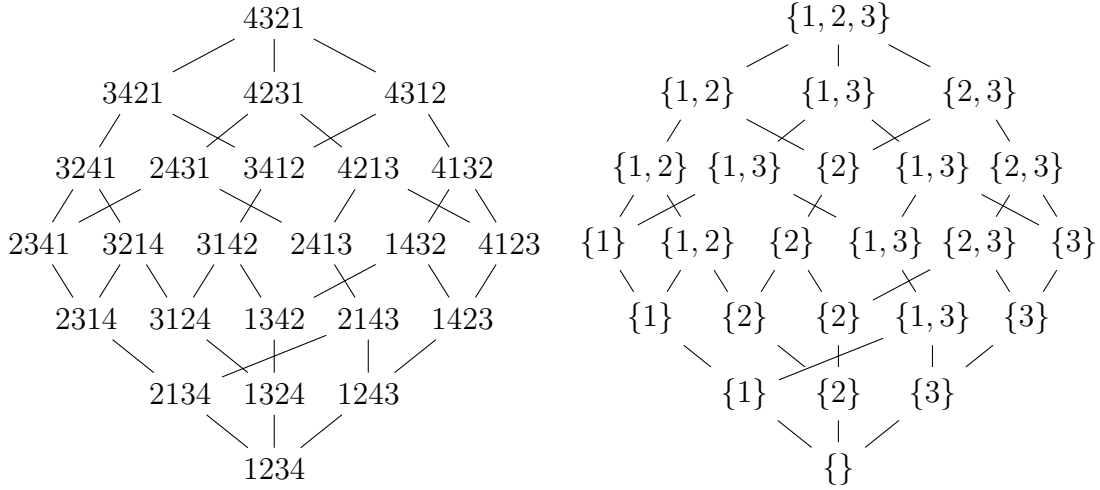

The (right) weak Bruhat ordering on the symmetric group produces a
lattice $\Bn = (\symn, \rel)$. For any permutation $p\in \symn$ we
write $p=p_1p_2\cdots p_n$ to mean $p(i)=p_i$ (this is the second line
in the two line notation).  The covering relation $\crel$ that generates the
relation $\rel$ is defined by $p_1p_2\cdots p_n\ \crel\ q_1q_2\cdots q_n$
whenever there is some $i \in \{1,\dots, n-1\}$, such that
$p_i<p_{i+1}$, $q_i=p_{i+1}$, $q_{i+1}=p_i$, and $q_k=p_k$
otherwise. This means that $p\ \crel\ q$ if $q$ is obtained by reversing
two consecutive and increasing elements from $p$. See
Figure~\ref{figS4} for $\Bfour$. 

The transpositions $(i,i+1)$ for $i \in \{1, \dots, n-1\}$ are called
the {\it generators}, and we will denote the generator $(i,i+1)$ by
$g_i$. We also define $\bG_n=\{g_i\mid 1\le i<n\}$ to be the {\it
  generating set} for $\Bn$. From another point of view, $\Bn$ is
the Cayley graph generated by $\bG_n$, turned into a lattice with
ranks given by distance from the identity $\bid$. The set of all
permutations at rank $\ell$ in $\Bn$ is denoted by  $\Bln$.

For a permutation $p\in\symn$ we define its set of {\it inverse
  descents} by
\[
  \rev(p)=\{ g_i \, \mid \,   i+1\textrm{ precedes }i\textrm{ in }p\},
\]
and for a set of permutations $\cP$, we define
$\rev(\cP)=\{\rev(p)\mid p\in\cP\}$. The set of inverse
descents play an important role in the Bruhat lattice, since for any $p \in \symn$,
we have $g_i\in \rev(p)$ if and only if $g_i\rel p$ in the weak order.

We denote the {\it inversion set} for a permutation $p \in \symn$ by
\[
\Inv(p)=\{ ( p_j,p_i )\, \mid \, 1\le i<j\le n, \, p_j<p_i\}
\]
with this definition, $\rank(p)=|\Inv(p)|$. Our definition of inverse
descents is nonstandard; usually the inverse descents are integers,
but we are using the simple transpositions indexed by those integers.

The inversion sets of different permutations are always distinct,
while their inverse descent sets are sometimes not. For example,
\begin{equation}\label{invertedgens}
\begin{aligned}
\rev (3214) &= \{g_1,g_2\}, \quad  & \Inv(3214) &= \{(1,2),(1,3),(2,3)\}, \\
\rev(3241) &= \{g_1,g_2\}, \quad &  \Inv(3241) &= \{(1,2),(1,3),(1,4),(2,3)\}.
\end{aligned}
\end{equation}

For a set $A \sse \bG_n$, define the {\it multiplicity of $A$} by
$| \{p\in \Bn \, | \, \rev(p) = A\}|$ and the {\it multiplicity of $A$
  at level $\ell$} by $| \{p\in \Bln \, | \, \rev(p) = A\}|$.  For
example, consider the set $\{g_1,g_4\}$. Its multiplicity in
$\Btwosix$ is $1$ since $21 3 54 6$ is the only permutation with rank
2 and inverse descents set exactly $\{g_1,g_4\}$. Its multiplicity in
$\Bthreesix$ is $3$ since each of $(1,3)$, $(3,5)$ and $(4,6)$ may be
reversed to increase the rank without introducing a new inverse
descent (the permutations in $\Bthreesix$ with inverse descents set
$\{g_1,g_4\}$ are $231 54 6$, $21 534 6$, and $21 3 564$).

\begin{lem}\label{lem:multi}
The multiplicity of any $k$-set of generators at level $\ell$ of $\Bn$
is at most $\binom{\ell+k-1}{k-1}e^{k\pi\sqrt{\frac{2\ell}{3}}}$.
\end{lem}
\proof Let $S=\{g_{i_1},g_{i_2},\ldots,g_{i_k}\}$ be a $k$-set of
generators and $\BlS$ be the set of permutations at level $\ell$ of
$\Bn$ whose inverse descent set is $S$. Set $i_0=0$, $i_{k+1}=n$, and
assume ${i_1<i_2<\cdots<i_k}$. Any permutation in $\BlS$ will preserve
the natural order of each of the sets
\[
S_0 = \{1,2,\dots, i_1\}, \, S_1:=\{i_1+1, \dots,  i_2\}, \dots , S_k=\{i_k+1,\dots,n\}
\]
and so we can identify every permutation in $\BlS$ with a permutation
of the multiset
\[
M_S=\{1^{i_1},2^{i_2-i_1},\ldots,k^{i_k-i_{k-1}},(k+1)^{n-i_k}\}.
\]
In general, there will be more multiset permutations than elements of
$\BlS$. An inversion in the multiset permutation represents an
inversion in the corresponding permutation in $\BlS$. The
$q$-multinomial coefficient
\begin{equation}
  \label{qmult}
  \qmult{n}{i_1}{i_2-i_1}{i_k-i_{k-1}}{n-i_k}=\qbin{n}{i_1}\qbin{n-i_1}{i_2-i_1}\cdots\qbin{n-i_k}{n-i_k}
\end{equation}
is the generating function by inversions for permutations of $M_S$ (see \cite{stanley}). Therefore, $|\BlS|$ is less than the coefficient of $q^{\ell}$ in \eqref{qmult}. The coefficient of $q^{\ell}$ in \eqref{qmult} is $\sum_{a_1+a_2+\cdots+a_k=\ell}\prod_{j=1}^kA_{a_j}$, where $A_{a_j}$ is the coefficient of $q^{a_j}$ in $\qbin{n-i_{j-1}}{i_j-i_{j-1}}$. In other words, 
\begin{equation}\label{bounda}|\BlS|\leq\sum_{a_1+a_2+\cdots+a_k=\ell}\prod_{j=1}^kA_{a_j}.\end{equation}.

Since $\qbin{n-i_{j-1}}{i_j-i_{j-1}}$ is the generating function by size for partitions whose Young diagram fits inside an $(i_j-i_{j-1})\times(n-i_j)$ rectangle (see \cite{stanley}), $A_{a_j}$ is the number of partitions of $a_j$ with at most $(i_j-i_{j-1})$ parts whose first part is at most $(n-i_j)$. We can crudely bound $A_{a_j}$ by $p(\ell)$, the number of partitions of $\ell$, and replace the bound in \eqref{bounda} with

\begin{equation}
  \label{boundb}
  |\BlS|\leq\sum_{a_1+a_2+\cdots+a_k=\ell}(p(\ell))^k.
  \end{equation}
The number of weak compositions of $\ell$ with $k$ parts is $\binom{\ell+k-1}{k-1}$ \cite{stanley} and $p(\ell)\leq e^{\pi\sqrt{\frac{2\ell}{3}}}$ \cite[Chapter 5]{andrews}, so we can replace \eqref{boundb} by

$$|\BlS|\leq \binom{\ell+k-1}{k-1}e^{k\pi\sqrt{\frac{2\ell}{3}}},$$

which is our result.
\pf

The bound above is often much larger than the actual multiplicity of a
set. The next result gives a much stronger bound for the number of
permutations in a level with the property that every inversion is a
generator. A set of inverse descents is called a \textsl{separated
  set} if it does not contain generators $g_i$ and $g_{i+1}$ for any
$i$.

\begin{prp}\label{prop:separated}
If a permutation $p$ in level $\ell$ of $\Bn$ has an
inverse descents set of size $\ell$, then the inverse descents set
must be a separated set.
\end{prp}
\proof
Assume that $A$ is a set of size $\ell$ that is not a separated set. 
Then $i,i+1 \in A$ for some $i$. Any permutation that has both $i$ and
$i+1$ in its inverse descents set must also have $(i,i+2)$ in its inversion
set. Thus the inversion set for any $p$ with $\rev(p) = A$ must have
at least $\ell +1$ elements, and $p$ is above rank $\ell$. So the
multiplicity of $A$ in $\Bln$ is 0.
\pf

The next claim follows from directly from the definitions of the
relation in the Bruhat order, inverse descents sets and inversion sets.

\begin{clm}\label{ordpres}
Let $p,q \in \Bn$, then the following hold.
\begin{enumerate}
\item If $p \rel q$, then $\rev(p) \sse \rev(q)$.\label{ordpres:four}
\item  $p \rel q$ if and only if $\Inv(p) \sse \Inv(q)$.\label{ordpres:one}
\item $p\wedge q\not=\bid$ if and only if $\rev(p)\cap\rev(q)\not=\mt$. \label{ordpres:three}
\end{enumerate}
\end{clm}

The converse of Statement~(\ref{ordpres:four}) in the above
claim is not true, since it is possible to have incomparable
permutations $p$ and $q$ with $\rev(p) = \rev(q)$. For example, $2431$
and $4213$ are incomparable, but have the same set of inverse
descents.

\begin{cor}\label{ordpres:two}
Let $p,q \in \Bn$, then $\rev(p\wedge q)=\rev(p)\cap\rev(q)$.
\end{cor}
\proof 
Statement~\ref{ordpres:four} implies that
$\rev(p\wedge q) \sse \rev(p)\cap\rev(q)$. Conversely, if $g_i\in\rev(p)\cap\rev(q)$, then $(i,i+1)\rel p\wedge q$, and
$g_i \in\rev(p\wedge q)$. Thus $\rev(p\wedge q)=\rev(p)\cap\rev(q)$
and the result follows. 
\pf

For a subset $A\sse \bG_n$, let
\[
\pi(A)=\bigvee_{{g_i}\in A}g_i
\]
this permutation is called the \textsl{minimal} element for $A$. Since
$\Bn$ is a lattice, $p$ is well-defined and unique.

\begin{lem}\label{witness}
  For any subset $A \sse \bG_n$ in the Bruhat ordering $p=\pi(A)$ is the
  unique minimum permutation such that $\rev(p) = A$. Further, for
  any $q$ with $\rev(q) =A$, both $\Inv(p) \sse \Inv(q)$ and $p \rel
  q$ hold.
\end{lem}
\proof Since $\pi(A)=\vee_{{g_i}\in A}g_i$, it follows that
$\rev(p) =A$.  Let $q$ be any other permutation with $\rev(q) =
A$. Then for any $g_1 \in A$ it follows that $g_1 \rel q$, and so
$p \rel q$ and, by Claim~\ref{ordpres} (\ref{ordpres:one}),
$\Inv(p) \sse \Inv(q)$.  \pf

Lemma~\ref{witness} implies that
\[
\pi(A)= \min\{p\in\Bn \, \mid \, \rev(p)=A\}.
\]
 
Lemma~\ref{witness} also shows that the multiplicity for any inverse
descents set is at least 1. However, it is possible to have a set that
is not an inversion set for any permutation; for example there is no
permutation with the inversion set $\{(1,2),(2,3)\}$ (if $2$ is before
$1$ and $3$ is before $2$, then $3$ is also before $1$).

The goal of the next section is to give properties of
intersecting families in the Bruhat lattice in order to
characterize the extremal intersecting sets.

%
%
\section{Properties of intersecting sets in $\Bn$}\label{BegRes}

The Bruhat ordering on permutations is not a uniquely
complemented lattice, but the permutations can be paired so that the
pairs have no common elements in either their inverse descents or
their inversion sets. The permutation $p=p_1p_2\cdots
p_n$ is paired with the permutation $\overline{p} = p_n\cdots p_2p_1$.
These pairs have the property that 
\[
\Inv(\op)=\{ (i,j)\, \mid \, 1\le i<j\le n\} \setminus \Inv(p),
\]
and also that $\rev(\op)=\bG_n\setminus\rev(p)$.  

With this paring, a version of Theorem~\ref{thm:complemented} holds for the Bruhat lattice.
We alter notation slightly to write $f_t(\Bn)$ as $f_t(n)$ and
$f_t(\Bkn)$ as $f_{t,k}(n)$.

\begin{thm}
Let $n$ be an integer with $n \geq 2$, then $f_1(n)=n!/2$.
\end{thm}
\proof Let $x$ be a permutation at level 1 in $\Bn$; assume that
$\rev(x) = \{g_i\}$. The upset of $x$ is the set of all permutations that
have $i$ in its inverse descents set. For each $p \in \Bn$, exactly
one of $p$ and $\op$ will have $i$ in its inverse descents set.  \pf

\begin{cor}
The Bruhat lattice has the 1-EKR property.
\end{cor}

The next goal is to give a characterization of the maximum
intersecting families in $\Bn$. If $\cP\sse\Bn$ is an intersecting
family of permutations of maximum size, then $\cP = n!/2$ and so, for
every permutation $p$, either $p \in \cP$ or $\op \in \cP$.

An intersecting family $\cP$ is called {\it maximal} if every permutation that intersects all the permutations in $\cP$ is also in $\cP$; equivalently, there are no elements not in $\cP$ that intersect all the elements in $\cP$.

For a family $\cP\sse\Bn$, define
\[
\min(\cP)=\{ p\in\cP \, \mid \,  \textrm{ if there exists $q \in \cP$ with $q\rel p$, then  $q=p$ } \}.
\]  
Observe that $\cP$ is an antichain if and only if $\min(\cP)=\cP$. Further, if
$\cP$ is a maximal intersecting family, then $\up(\min(\cP))=\cP$. 

\begin{lem}\label{lem:maxinter}
  Let $\cP\sse\Bn$ be a maximal intersecting family.  If $p \in
  \min(\cP)$, then $p = \pi(\rev(p))$. 
\end{lem}
\proof
Let $p \in \min(\cP)$ and set $A = \rev(p)$. Then for any
permutation $q \in \cP$ we have that $\rev(q) \cap A \neq \emptyset$.
Thus $\pi(A)$ intersects every $q \in \cP$. Since $\cP$ is maximal, $\pi(A) \in \cP$.
By definition of $\pi(A)$ we have $\pi(A) \rel p$ and, since $p$ is
minimal in $\cP$, we have $p = \pi(A)$.
\pf

In developing our characterization of maximum intersecting families we will move between permutations and their sets of inverse descents, to do this we will use the following lemma.

\begin{lem}\label{LemEquiv}
Let $\cP\sse\Bn$. Then, the following statements are equivalent.
\begin{enumerate}
\item
$\cP$ is intersecting.\label{One}
\item
$\up(\cP)$ is intersecting.\label{Two}
\item
$\min(\cP)$ is intersecting.\label{Three}
\item
$\rev(\cP)$ is an intersecting set system.\label{Four}
\item 
If $\cG$ is a system of subsets of $\bG_n$, with
  $\cP=\{\pi(g) \, | \, g \in \cG \}$, then $\cG$ is intersecting.\label{Five}
\end{enumerate}
\end{lem}

\proof Since $\Bn$ is a poset, for all $p,p^\pr, q, q^\pr\in\cP$
with $p\rel p^\pr$ and $q\rel q^\pr$, it follows that $p\wedge q\rel
p^\pr\wedge q^\pr$. This implies that if a family $\cP$ is intersecting,
then $\up(\cP)$ is also intersecting.  This fact shows that Statement
(\ref{One}) implies Statement (\ref{Two}).  Statement (\ref{Two})
implies Statement (\ref{Three}) since $\min(\cP)
\sse\up(\cP)$. Statement (\ref{Three}) implies Statement (\ref{One})
since $\cP \sse \up(\min(\cP))$. Thus the first three statements are
equivalent.

Finally, Statement~\ref{ordpres:three} of Claim~\ref{ordpres} implies that Statements (\ref{Four}) and
(\ref{Five}) are equivalent to (\ref{One}).  \pf

For any set system $\cA$ defined on an $(n-1)$-set, which we may take to be $\bG_n$,
there is a corresponding subfamily of $\Bn$ defined by
\[
\sP(\cA)=\up \left(\{\pi(A) \mid A \in \cA\} \right).
\]
This means that $\sP(\cA)$ is the family of all permutations $p$ such
that there exists an $A \in \cA$ such that $A \sse \rev(p)$.

\begin{cor}\label{cor:SetsImpliePerms}
  If $\cA$ is an intersecting set system on $\bG_n$, then
  $\sP (\cA)$ is an intersecting family in the Bruhat ordering.
\end{cor} 
\proof
This follows from Statement~\ref{ordpres:three} of Claim~\ref{ordpres} and Part~(\ref{Five}) of Lemma~\ref{LemEquiv}.
\pf

Further, it is possible to construct a set system from a family of
permutations. If $\cP$ is a family of permutations, then
$\rev(\min(\cP))$ is a set system on $\bG_n$ called the
\textsl{generating set for $\cP$.} The next three lemmas give
properties of this set system.

\begin{lem}\label{MaxInt}
  Let $\cP\sse\Bn$ be an intersecting family and define $\cA
  = \rev(\min(\cP))$.  Then $\cA$ is
  an intersecting set system.  Further, if $\cP$ is maximal then  $p \in \cP$ if and only if $A \sse
  \rev(p)$ for some $A\in\cA$.
\end{lem}

\proof The first statement is clear from Statement~\ref{ordpres:three} of Claim~\ref{ordpres}.

The implication in the second statement follows from the definition of
$\cA$. In particular, if $ p \in \cP$ then $A \sse \rev(p)$ for some
$A\in\cA$. Conversely, for any permutation $p$, if $A \sse \rev(p)$
for some $A\in\cA$ then $\pi(A) \in \cP$ and $\pi(A) \rel p$. Since
$\cP$ is maximal, this implies that $p \in \cP$.  \pf

\begin{cor}\label{cor:maxP}
Let $\cP\sse\Bn$, then $\cP \sse \sP( \rev(\min(\cP)) )$. If $\cP$ is a maximal intersecting
family, then $\cP = \sP(\rev(\min(\cP)))$.
\end{cor}

The next result shows that $\rev$ is a weak order
preserving map on the poset, so it preserves the $\min$ operation.

\begin{cor}\label{mininv}
  Let $\cP\sse\Bn$ be a maximal intersecting family of permutations, then
  $\rev(\min(\cP))=\min(\rev(\cP))$.
\end{cor}
\proof From Part~(\ref{ordpres:four}) of Claim~\ref{ordpres}, if $p
\rel q$, then $\rev(p) \sse \rev(q)$. This implies that if $\rev(p)
\in \min(\rev(\cP))$, then $p \in \min(\cP)$. Thus $\min(\rev(\cP))
\sse \rev(\min(\cP))$.

Conversely, if $\rev(p) \in  \rev(\min(\cP))$, then $\rev(p) =
\rev(p')$ for some $p' \in \min(\cP)$. Since $\cP$ is maximal,
Lemma~\ref{lem:maxinter} implies that $p' = \pi(\rev(p'))$. 

If there is a set $\rev(q) \in \rev(\cP)$ with $\rev(q) \subset
\rev(p')$, then $\pi(\rev(q)) \rel \pi(\rev(p'))$, which contradicts $p' \in \min(\cP)$.
This means that $\rev(p')$ is a minimal set in $\rev(\cP)$ and 
thus $\rev(p) = \rev(p') \in \min(\rev(\cP))$.
\pf

Corollary~\ref{mininv} implies that $\rev(\min(\cP))$ is an antichain.
We call $\cA$ a \textit{maximal intersecting antichain} if $\cA$ is an
intersecting antichain with the additional property that, for any set
$B \not\in \cA$, there is a set $A\in\cA$ such that either
$A\cap B=\mt$, or $A\sse B$. Equivalently, an intersecting antichain
$\cA$ is a maximal intersecting antichain if, for any $B$ with
$X \cap B \neq \emptyset$ for all $X \in \cA$, there exists an
$A \in \cA$ with $A \subseteq B$.

\begin{lem}\label{lem:MaxAntiChains}
Let $\cP\sse\Bn$ be a maximal intersecting family and define $\cA = \rev(\min(\cP))$.
Then $\cA$ is a maximal intersecting antichain.
\end{lem}
\proof From Corollary~\ref{MaxInt}, $\cA$ must be intersecting.
Since $\cA=\rev(\min(\cP))$, Corollary~\ref{mininv} implies $\cA=
\min(\rev(\cP))$ and that $\cA$ is
an antichain.

Suppose that $\cA$ is not a maximal intersecting antichain. Then there
must be some $B\not\in\cA$ such that $\cA \cup \{B\}$ is intersecting
and there is no $A \in \cA$ with $A \sse B$.  Let $p = \pi(B)$.  Since
no $A \in \cA$ has the property that $A \sse \rev(p)$, we have that
$p \not \in \cP$. Then $\cP \cup \{p\}$ is a intersecting family,
contradicting the maximality of $\cP$.\pf

At this point, to avoid unnecessary technicalities, we associate any set of 
integers with the set of generators having those integers as subscripts;
e.g. $\{1,3,4\}\equiv\{g_1,g_3,g_4\}$.

There are many maximal intersecting set systems that are not maximal intersecting
antichains. For example, consider the following system of
$r$-subsets of $M=\{1,\ldots,m\}$ (with $m >2r$) defined by
\begin{eqnarray}\label{eq:Hnrk}
\cH_{m,r,k}=\{H\in\binom{M}{r}\mid |H\cap \{1,\dots, 2k-1\} | \ge k \},
\end{eqnarray}
where $1\le k \le r\le (m-1)/2$. If $k=1$, this is the canonical
intersecting set. But the set system $\cH_{m,r,k}$ 
is a maximally intersecting antichain if and only if $k=r$. If $r>k$, consider the set
$\{r,\ldots,m\}$, which intersects each set in $\cH_{m,r,k}$. Since $r>k$,
the set $\{r,\ldots,m\}$ is neither contained in, nor contains any of the
sets from $\cH_{m,r,k}$. The following sets are examples of maximally intersecting antichains:
\begin{align*}
\cA = & \{ A \in \binom{M}{r} \, \mid \, |A \cap \{1,2,...,2r-1\}| \geq r\}, \\
\cB = &  \{12,\; 13,\; 14, \;15, \cdots, 1m\;, 234\cdots m \}, \textrm{\ and}\\
\cC = &  \{123, 134, 145,156,126,235,245,246,346,356\} \textrm{ on an $6$-set.}
\end{align*}
The sets below are not maximally intersecting antichains:
\begin{align}
\mathcal{D} = & \{ A \in \binom{M}{r} \,\mid \, |A \cap \{1,2,...,2k-1\}| \geq k\}, \quad \textrm{if }  r > k \nonumber, {\rm\ and} \\
\cE = & \{ A \in\binom{M}{r} \, \mid \, 1 \in A \}, \quad \textrm{ if } r>1 \nonumber.
\end{align}

The final theorem of this section is our characterization of the
maximum intersecting families in the Bruhat lattice.

\begin{thm}\label{MIGthm}
Let $\cP\sse\Bn$ be an intersecting family and define $\cA = \rev(\min(\cP))$. Then $\cP$ is maximum if and only if $\cA$ is a maximal intersecting
antichain and $\mathcal{P} = \mathcal{P}(\mathcal{A})$.
\end{thm}

\proof We have already seen in Lemma~\ref{lem:MaxAntiChains} that if
$\cP$ is maximum, then $\cA$ is a maximal intersecting
antichain. Further, by Corollary~\ref{cor:maxP},
$\mathcal{P} = \mathcal{P}(\cA)$.

Assume that $\cA$ is a maximal intersecting antichain and
$\mathcal{P} = \mathcal{P}(\mathcal{A})$. We will show that for every
permutation $p$, either $p \in \cP$ or $\op \in \cP$.

For a permutation $p \in \Bn$ let $B = \rev(p)$.  If $B$ intersects
every set in $\cA$ then, by maximality, $A \sse B$ for some
$A \in\cA$.  This implies that $p \in \up(\pi(\cA)) \sse \cP$.  If $B$
does not intersect every set in $\cA$, then $B \cap A = \emptyset$ for
some $A\in\cA$. This implies that $A \sse \oB$ so
$\op\in \up(\pi(A)) \subseteq \cP(\mathcal{A}) = \cP$.  \pf

From Corollary~\ref{cor:SetsImpliePerms}, the family
$\sP(\cH_{m,r,k})$ (defined in Equation~\eqref{eq:Hnrk}) is
intersecting.  There are many different set systems that are
isomorphic to $\cHmrk$, formed by permuting the underlying set.  It is
possible that isomorphic maximal intersecting antichains give rise to
non-isomorphic maximum intersecting families of permutations.  For
example, when $(m,r,k)=(6,2,2)$, then
$\cH_{6,2,2}=\{ (1,2), (1,3), (2,3)\}$, and
$\cH = \{(1,3),(1,5),(3,5)\}$ is isomorphic to $\cH_{6,2,2}$. Both set
systems can be used to construct an intersecting set of permutations
in the Bruhat order
\begin{align*}
\cP(\cH_{6,2,2}) &= \up(\{321456,  214356, 143256\}), \\
\cP(\cH)  &=  \up(\{214356, 213465,  124365\}).
\end{align*}
Both of these sets have size $360$, but they are non-isomorphic
families in $\Bn$. This can be seen since the ranks of
$321456$ and $143256$ are $3$, whereas the other permutations are of rank $2$.
One can see that the key distinction here is that the sets in $\cH$ are separated sets.

The final Corollary of this section follows from Theorem~\ref{MIGthm}.

\begin{cor}
  Let $\cH$ be a set system isomorphic to $\cH_{m,r,k}$, with $1\le k\le r<(m-1)/2$.
  If $k=r$ then $\cP(\cH) =m!/2$ while, if $k<r$, then $\cP(\cH) <m!/2$.
\end{cor}

%
%
\section{Erd\H{o}s--Ko--Rado Theorem for levels in the Bruhat lattice}

In this section we prove that an Erd\H{o}s--Ko--Rado Theorem holds
for the $r^\textrm{th}$ level of $\Bn$, provided that $n$ is large relative to $r$.
Recall that $\Brn$ denotes the set of all permutations of rank $r$ in the Bruhat lattice $\Bn$, and that $\Brn$ is EKR if the size of its largest intersecting subfamily is no larger than the largest 1-star at level $r$. It is clear that $\Btwon$ is EKR. In this section we first give a proof that $\Bthreen$ is EKR, and more generally that $\Brn$ is EKR, provided that $n$ is large relative to $r$.

\begin{lem}\label{Rank3StarSize}
Let $\cH=\Bthreen$ and $p\in\Bonen$. Then $|\up_\cH(p)|=\binom{n-1}{2}$.
\end{lem}

\proof 
Define the polynomial $[m]=(1+x+\cdots+x^{m-1})$; then the generating function for the sizes of the ranks in $\Bn$ is
\[
F(x)=[n]!=(1+x)(1+x+x^2)\cdots (1+x+x^2+\cdots +x^{n-1}).
\]
From~\cite{wei}, the generating function for the sizes of the ranks in $\up(p)$, where $\rank(p)=1$, is given by the generating function
\[
\frac{F(x)}{1+x}=\frac{[n]!}{1+x}=(1+x+x^2)\cdots (1+x+x^2+\cdots
+x^{n-1}).
\]
Hence $|\up_\cH(p)|$ is the coefficient of $x^2$ in $F(x)/(1+x)$. This
coefficient is the same as that in $(1+x+x^2)^{n-2}$, which is
$\binom{n-2}{2}+(n-2)=\binom{n-1}{2}$.  \pf

We will use the following result of Holroyd, Spencer and
Talbot~\cite{HST} on EKR theorems for separated sets. 

\begin{thm}[\cite{HST}]\label{HSTthm}
  The largest intersecting system of separated $r$-sets from an $m$-set
  has size $\binom{m-r}{r-1}$ and consists of all the separated sets
  that contain some fixed point. \pf
\end{thm}

\begin{thm}\label{Rank3EKR}
The set $\Bthreen$ is EKR for all $n\ge 3$.
\end{thm}
\proof 
If $n=3$ the result is evident, so we let $n\ge 4$.
Let $\cP\sse\Bthreen$ be an intersecting family of maximum
size. We are done if $\cP$ is a star, so we assume otherwise.
Let $\cP_i=\{p\in \cP \, | \, |\rev(p)|=i\}$ where
$i \in \{1,2,3\}$. Clearly, $\cP=\bigcup_{i=1}^3 \cP_i$.

If $p \in \cP_1$, then, since $\cP$ is intersecting, every permutation in $\cP$ has the 
single element in $\rev(p)$ in its inverse descent set. This implies that $\cP$ is a star and we are done. So we
may assume that $\cP_1 = \emptyset$. Now we will consider $\cP_i$ for
$i >1$.

Any permutation $p$ in $\cP_3$ must have $\rev(p) = \Inv(p)$, since it
is a rank 3 permutation with 3 inverse descents, and by Proposition~\ref{prop:separated}
$\rev(p)$ must be a separated set and every set $\Inv(p)$, for
$p \in \cP_3$, has multiplicity $1$.  Thus the set system $\{\rev(p)
\, | \, p \in \cP_3\}$ is an intersecting system of separated sets
that has the same size as $\cP_3$.  By Theorem~\ref{HSTthm} (with
$m=n-1$) we have $|\cP_3|\le \binom{n-4}{2}$.

Next consider $\cP_2$ and define the set system
$\cA_2 = \{\rev(p)\ |\ p \in P_2 \}$. Clearly, $\cA_2$ is an
intersecting $2$-set system. If $|\cA_2| \geq 4$, then $\cA_2$ must
contain four sets isomorphic to
\[
\{12,13,14,15\}.
\]
Since the inverse descent sets for the permutations in $\cP_2$ and
$\cP_3$ are sets of size no more than 3, and these inverse descents
sets intersect all the sets in $\cA_2$, this implies that all the
permutations in both $\cP_2$ and $\cP_3$ must have a common element in
their inverse descent. Thus if $|\cA_2| \geq 4$, then $\cP$ is a
star. So we can assume that $|\cA_2|\leq 3$. The multiplicity of any
set in $\cA_2$ is no more than $4$, so $|P_2| \leq 12$.

Since $\cP$ is not a star, we have
\begin{equation}\label{eq:bound3}
|\cP| \leq 12 + \binom{n-4}{2} = \frac{1}{2} (n^2-9n+44),
\end{equation}
which is no more than
$ \binom{n-1}{2}  = \frac{1}{2}(n^2-3n+2)$,
provided that $n \geq 7$.

For $n\leq 6$, Equation~\eqref{eq:bound3} provides a bound which
is at most the size of a star in the poset. If $n=4$ then, from
Figure~\ref{figS4}, each set in $\cA_2$ has multiplicity $2$. The
bound in Equation~\eqref{eq:bound3} becomes $|\cP| \leq 3$, which is
exactly the size of a star in the poset. Similarly, when $n=5$ the
multiplicities of the sets in $\cA_2$ are no more than $2$ and the
bound becomes $|\cP| \leq 6$. This again is the size of a star in the
poset. Finally, when $n=6$, the multiplicity is no more than 3 and the
bound is $|\cP| \leq (3)(3)+1 = 10$. This also is the size of a star
in the poset.  \pf 

Next we will use a counting method, similar to the one used in
Theorem~\ref{Rank3EKR}, to show that $\Brn$ is EKR, provided that $n$
is sufficiently large. First, we will prove a lower bound on the size
of a star in $\Brn$, and then give two known results that can be used
to show that if an intersecting set is sufficiently large, then there
are minimum number of sets that intersect in exactly one element.

\begin{lem}\label{RankrStarSize}
Let $\cH=\Brn$ and $p\in\Bonen$. Then $|\up_\cH(p)| > \binom{n-2}{r-1}$.
\end{lem}
\proof
As we have seen, if $p \in \Bonen$, then the size of $\up_\cH(p)$ is given by the coefficient of $x^{(r-1)}$ in 
\[
\frac{[n]!}{x+1} = (1+x+x^2)(1+x+x^2+x^3) \cdots (1+x+x^2+\cdots +x^{n-1}).
\]
There are $\binom{n-2}{r-1}$ ways to select the term $x$ from $r-1$
different factors. So the coefficient of $x^{(r-1)}$ is at least
$\binom{n-2}{r-1}$.  \pf

The next result is the Hilton-Milner theorem. This result gives a
bound on the size of an intersecting set system in which the sets do
not all contain a common element.

\begin{thm}[\cite{MR0219428}]\label{thm:hm}
If $\mathcal{A}$ is an intersecting system of $k$-subsets of an $m$-set with $\cap_{A \in \cA} A = \emptyset$, then
\[
|\mathcal{A} | \leq \binom{m-1}{k-1} - \binom{m-k-1}{k-1} + 1.
\]
\end{thm}

We will also use a version of Erd\H{o}s' matching conjecture due to Frankl.

\begin{thm}[\cite{MR3033661}]\label{thm:frankl}
  Let $\mathcal{A}$ be a system of $k$-subsets of an $m$-set. Assume
  that $\mathcal{A}$ contains no $r+1$ pairwise disjoint sets and that
  $m \geq (2r+1)k-r$. Then
\[
|\mathcal{A}| \leq \binom{m}{k} - \binom{m-r}{k}.
\]
\end{thm}

The previous two results can be used to show that if an intersecting
set system is sufficiently large, then there will be many sets that
intersect in exactly in one common element.

\begin{lem}\label{lem:bigIntersectingFams}
  Let $\mathcal{F}$ be an intersecting system of
  $\ell$-sets on an $(n-1)$-set with $2 \leq \ell$ and $n \geq (2r+1)(\ell-1)-(r-2)$. If $\ell < r $ and
  \[
  |\mathcal{F} | > \binom{n-2}{\ell -1} -\binom{n-2-r}{\ell -1},
  \]
  then every set in $\mathcal{F}$ includes a fixed element $x$, and
  $\mathcal{F}$ contains at least $r+1$ sets that pairwise intersect
  in exactly the element $x$.
\end{lem}
\proof 
First note that since $2 \leq \ell < r$, it follows that 
\[
|\mathcal{F}| > \binom{(n-1) -1}{\ell-1} - \binom{(n-1) - (\ell+1)}{\ell-1} +1,
\]
and Theorem~\ref{thm:hm} implies that every set in $\mathcal{F}$
contains some element $x$. Every set in $\mathcal{F}$ contains $x$, so
let $\mathcal{F}^x$ be the family of $(\ell-1)$-sets on an $(n-2)$-set
formed by removing $x$ from every set in $\mathcal{F}$.  Applying
Theorem~\ref{thm:frankl} to $\mathcal{F}^x$ yields that, for the
stated values of $n$, the family $\mathcal{F}^x$ contains at least
$r+1$ sets that are pairwise disjoint.  \pf

The previous lemma will be used to bound the number of permutations
$p$, in an intersecting family in $\Brn$, that have $|\ID(p)|=\ell<r$.
Theorem~\ref{HSTthm} will be used to bound the number of such $p$
having $|\ID(p)|=r$, since in this case
Proposition~\ref{prop:separated} implies that $\ID(p)$ is a separated
set.

\begin{thm}\label{RankrEKR}
The set $\Brn$ is EKR for $n$ sufficiently large relative to $r$.
\end{thm}
\proof
We may assume that $r>3$, since we have seen that both $\Btwon$ and $\Bthreen$ are EKR.
Let $\cP\sse\Brn$ be an intersecting family of maximum size.
We are done if $\cP$ is a star, so we assume otherwise.
As in the proof for $\Bthreen$, let $\cP_i=\{p\in \cP \, | \, |\rev(p)|=i\}$ where
$i \in \{1,2,\dots,r\}$. Clearly, $\cP=\bigcup_{i=1}^r \cP_i$. Further set $\cA_i = \{\rev(p) \,|\, p \in P_i \}$.

As in the proof of Theorem~\ref{Rank3EKR}, if $\cP_1 \neq \emptyset$ then then $\cP$ is a star,
so we may assume that $\cP_1 = \emptyset$. Now we will consider
$\cP_\ell$ for $\ell \in \{2,\dots, r-1\}$.

By Lemma~\ref{lem:bigIntersectingFams}, if
$|\cA_\ell| >\binom{n-2}{\ell -1} -\binom{n-2-r}{\ell -1}$ then
$\cA_\ell$ is a star, centered on some $x$, which contains at least
$r+1$ sets that pairwise intersect at exactly $x$. Every $p\in\cP$ has
$|\ID(p)|\le r$ and intersects each of these $r+1$ sets. Thus each
$\ID(p)$ contains $x$, meaning that $\cP$ is a star. Hence, we may
assume that $|\cA_\ell| \leq \binom{n-2}{\ell -1} -\binom{n-2-r}{\ell -1}$ for
each $\ell \in \{2, \dots, r-1\}$. Finally, since $\cA_{r}$ is an intersecting
family of separated $r$-sets from an $(n-1)$-set, Theorem~\ref{HSTthm} implies that $|\cA_r| \leq \binom{n-1-r}{r-1}$.

The multiplicity of any set from $\cA_\ell$ in $\cP_\ell$, is at most $\binom{\ell+r-1}{\ell-1} e^{\pi\ell \sqrt{\frac{2r}{3}}}$. Thus
\begin{align} \nonumber
|\cP| &\leq \binom{n-1-r}{r-1}+ \sum_{\ell = 2}^{r-1}
    \binom{\ell+r-1}{\ell-1} e^{\pi\ell \sqrt{\frac{2r}{3}}} \left(\binom{n-2}{\ell -1} -\binom{n-2-r}{\ell -1} \right) \\
&\leq \binom{n-1-r}{r-1}+ \sum_{\ell = 2}^{r-1}\binom{\ell+r-1}{\ell-1} e^{\pi\ell \sqrt{\frac{2r}{3}}}
\left(\sum_{j=3}^{r+2} \binom{n-j}{\ell-2}\right) \label{eq:pascal} \\
& \leq  \binom{n-3}{r-1} + (r-2)\binom{2r-2}{r-2} e^{3r^{\frac{3}{2}}}(r) \binom{n-3}{r-3}.  \label{eq:bounds}
\end{align}

Equation~(\ref{eq:pascal}) follows from repeated applications of
Pascal's Rule. The bounds used in Equation~(\ref{eq:bounds}) are very
rough upper bounds.

Provided that $n$ is larger than
\[
(r-2) \binom{2r-2}{r-2} e^{3r^{\frac{3}{2}}}r(r-2) + r,
\]
the size of $\cP$ is less than $\binom{n-2}{r-1}$, which by Lemma~\ref{RankrStarSize} is less than the size of the 
largest star.
\pf

The previous result uses a very rough counting method to bound the
size of an intersecting family of permutations that is not a star. We
believe that the lower bound on $n$ given above is far larger than
necessary; in fact we conjecture that all levels below the middle
level in the weak Bruhat lattice have the EKR property.

\begin{cnj}
  For all values of $1\le r\le \frac12 \binom{n}{2}$ the set $\Brn$ is
  EKR.  Thus, for any $p\in\Bonen$ we have $f_{1,r}(n)=|\up(p)|$,
  which equals the coefficient of $x^{r-1}$ in $F(x)/(1+x)$.
\end{cnj}

Figure~\ref{figS4} reveals that the conjecture holds for $\Bfour$. The
upper bound on $r$ is based on the fact that permutations in levels
above $\frac12 \binom{n}{2}$ have larger sets of inverse descents. For
sufficiently large $r$, the entire $r$-level is intersecting since the
size of the inverse descent sets are large. We suspect that for levels
above $\frac12 \binom{n}{2}$ the sets of all permutations with inverse
descent sets of size at least $(n-1)/2$ are larger than a star.

\section{$t$-intersecting families}\label{tIntFam}

In this section we conjecture that a version of the EKR theorem holds
for some levels of the Bruhat lattice for $t$-intersecting permutations. The original
version of the EKR theorem was for $t$-intersecting sets, and can be
stated as follows.

\begin{thm}\label{thm:ekrt}
  Let $k,t$ and $n$ be positive integers with $t<k$. Assume that $\cA$
  is a $t$-intersecting system of $k$-subsets of $\{1,\dots,n\}$.
  There exists a function $f(k,t)$ such that, if $n >f(k,t)$,
  then
	\[
		|\cA| \le \binom{n-t}{k-t}.
	\]
	Moreover, equality holds if and only if $\cA$ is a $t$-star at level $k$.
\end{thm}  

Similar to the EKR theorem for $t$-intersecting $k$-sets, we conjecture
that the largest $t$-intersecting family of rank $k$ permutations
in the Bruhat lattice is the upset of a permutation at level
$t$. Unlike the case for sets, the sizes of the upsets of the
permutations at level $t$ are not all equal in the Bruhat lattice. Before stating the
conjecture, we need to determine which permutation at level $t$ we
believe has the maximum number of rank $r$ permutations in its upset.

Consider the permutation $\rho(t)$ that is formed by starting with the trivial permutation
$\langle 1,2,\dots, n\rangle$ and reversing exactly $t$ pairs in $\{1,\ldots,n\}$. These reversals are done
iteratively in $t$ steps. At each step, the largest element in
the permutation that has a smaller element directly before it (this is
the largest element with a smaller element on its left), is switched with the
element directly to the left of it. 
For example, for $n=6$ we have
\begin{align*}
\rho(0) &= 123456, &\rho(1) &= 123465, &\rho(2) &= 123645, &\rho(3) &= 126345,  \\
\rho(4) &= 162345, &\rho(5) &= 612345, &\rho(6) &= 612354, &\rho(7) &= 612534, \\
\rho(8) &= 615234, &\rho(9) &= 651234, &\rho(10) &= 651243, &\rho(11) &= 651423, \\
\rho(12) &= 654123, &\rho(13) &= 654132, &\rho(14) &= 654312, &\rho(15) &= 654321.
\end{align*}

If $i$ and $j$ are positive integers with $i$ chosen as large as
possible so that
\begin{equation}\label{eq:defij}
t = (n-1) + (n-2) + \dots + (n-i) + (n-(i+j+1)) 
\end{equation}
then
\[
\rho(t) = \langle n, n-1, \dots, n-i+1, 1, 2, \dots, j, n-i, j+1, j+2, \dots, n-i-1\rangle.
\]
The set
$\Inv(\rho(t))$ is the set of the final $t$ transpositions in the
lexicographic ordering 
(for example, if $n=6$, then $\Inv(\rho(11))=\{(2,4),(3,4),(1,5),\ldots,(4,5),(1,6),\ldots,(5,6)\}$)
and $\rho(t)$ is at level $t$ in the Bruhat lattice. Further, the
inverse descent set for $\rho(t)$ is the set of the final $i+1$ generators.

By Corollary~3.7 of~\cite{wei}, the generating function for the sizes
of the levels in the down set of $\rho(t)$ is given by
\[
[n][n-1] \cdots [n-i+1] [n-i-j],
\]
(where $i$ and $j$ are as defined in Equation~\ref{eq:defij}).
By Theorem 3.1 of~\cite{wei} this implies that generating function of
the levels in the upset of $\rho(t)$ is
\[
F_{n,t}(x)\ :=\ [n-i] \cdots [n-i-j+1][n-i-j-2] \cdots [2][1].
\]
For example, with $n=4$ and $t=4$ (see Figure \ref{figS4}), we have
$\rho(4)=4132$, $i=1$, and $j=1$.  This yields the generating
functions $[4][2]=1+2x+2x^2+2x^3+x^4$ for its downset and
$F_{4,4}(x)=[3]=1+x+x^2$ for its upset.

The inverse descent set of $\rho(t)$ is very small. This means that
there are many different permutations of a given rank with an inverse
descent set that contains $\Inv(\rho(t))$. We conjecture that
$\rho(t)$ is a rank $t$ permutation with the maximum number of rank
$r$ permutations in its upset.

\begin{cnj}
  Let $1\le t\le r\le \binom{n}{2}$. If $n$ is sufficiently large
  relative to $t$ and $r$, then the set $\Brn$ is $t$-EKR.  In
  particular, $f_{t,r}(n)=|\up(p) \cap \Brn|$, which equals the
  coefficient of $x^{r-t}$ in $F_{n,t}(x)$.
\end{cnj}

Finally we note that the question of $t$-intersection for any $t<n$
can also be considered for the entire poset. In the case of the
subsets poset, the largest $t$-intersecting sets are given by Katona's
theorem~\cite{MR0168468}.  The theorem states that the largest such
set is roughly the collection of all set with size greater than
$(n+t)/2$. At present it is not clear to us what the case will be for
the Bruhat poset. Two possible candidates for largest $t$-intersecting
family are the set of all permutations with at least $(n+t)/2$ inverse
descents, and $\up(\rho(t))$ (this is the set of all permutations that have a common set of
size $t$ in their inverse descent set).

\begin{qst}
What is the size and structure of the largest $t$-intersecting
permutations in $\Bn$?
\end{qst}

\section*{Acknowledgment}
Karen Meagher's research is supported in part by NSERC Discovery
Research Grant Application number RGPIN-341214-2013.  Susanna Fishel's
is partially supported by a collaboration grant for mathematicians
from the Simon Foundation, grant number 359602.  Glenn Hurlbert's
research is supported in part by the Simons Foundation, grant number
246436. Glenn Hurlbert and Vikram Kamat wish to acknowledge the
support of travel grant number MATRP14-17 from the Research Fund
Committee of the University of Malta.  This grant enabled them to
travel to the 2MCGTC2017 conference in Malta where some of this
research was completed. The authors also wish to thank the referees
for useful comments that helped improve this article.

%
%

\[\]
%
 \end{document}